%% file: paper.tex
\documentclass[10pt,reqno]{amsart}

\input{research-notes-preamble.tex}

\usepackage{bm}

\DeclareMathOperator{\Br}{Br}

\newcommand{\para}[1]{\paragraph{\color{Black}\textbf{#1}}}

\usepackage[mode=tex]{standalone}% requires -shell-escape

\newif\ifhideexers
\ifhideexers
\usepackage{environ}
\NewEnviron{hide}{}

\fi

\graphicspath{ {resources} }

\newcommand{\VARTitle}{Torsion in the Braid Monodromy of Elliptic Fibrations} 
\newcommand{\VARAuthor}{Faye Jackson}
\author{\vspace{-0.25cm}\VARAuthor}

\title{\vspace*{-2cm}\VARTitle} % Change title as needed
\date{}
\definecolor{darkpastelpurple}{rgb}{0.59, 0.44, 0.84}
\definecolor{mauve}{rgb}{0.88, 0.69, 1.0}

\begin{document}

\hypersetup{linkcolor=lapislazuli}

\begin{abstract} 
	Given an elliptic fibration $\pi : M \to S^2$ with singular locus $\Delta \subseteq S^2$, let $\Br(\pi) < \Mod(S^2,\Delta)$ be the subgroup of the spherical braid group consisting of those braids that lift to a fiber-preserving diffeomorphism of $M$. We classify the order $n = \abs{\Delta}$ elements of $\Br(\pi)$ up to conjugacy in $\Br(\pi)$. To do so, we relate these conjugacy classes to special points on the $\SL_2$-character variety for $(S^2,\Delta)$ that correspond naturally to the exceptional elliptic curves $\CC/\ZZ[\omega]$ and $\CC/\ZZ[i]$ with their associated norms on $\ZZ[\omega]$ and $\ZZ[i]$. We also show that there are no elements of order $n-1$ or $n-2$ in $\Br(\pi)$, as there are in $\Mod(S^2,\Delta)$.
\end{abstract}

\thispagestyle{empty}
\maketitle
\setcounter{tocdepth}{2}
\microtypesetup{protrusion=false}
%\tableofcontents
\microtypesetup{protrusion=true}
%\newpage

\lhead{\VARAuthor}
\chead{\VARTitle}
\rhead{}
\setstretch{1.2}

\section{Introduction}

\subfile{intro.tex}

\para{Acknowledgements}

I thank my advisor Benson Farb for his substantial effort in improving this paper and his support throughout this project. Jointly, I thank Benson Farb and Eduard Looijenga for spending time explaining their work to me. I am particularly grateful to Seraphina Lee, for her excellent thesis presentation which inspired this paper. Thank you also to Dan Margalit, for his comments on an initial draft on this paper which made the introduction clearer. Finally, I thank the members of the geometry and topology working group for their time and interest in my project, particularly Ethan Pesikoff and Zhong Zhang. 

\section{Monodromy and Lifting Mapping Classes}\label{sec:prelim}

Let $\pi : M \to S^2$ be an elliptic fibration. We assume throughout that each singular fiber of $\pi$ has one critical point i.e., that the singularities are stable. For convenience, let $\Delta_\pi \subseteq S^2$ denote the set of singular values of $\pi$;  when the fibration $\pi$ is unambiguous we will instead write $\Delta = \Delta_\pi$. As discussed in the introduction, the braids in $\Mod(S^2,\Delta)$ which lift to a fiber-preserving diffeomorphism of $M$ are referred to as the \textit{liftable braids}, and the subgroup consisting of liftable braids is denoted $\Br(\pi)$. 

We now give a description of the liftable braids in terms of the Hurwitz action, essentially due to Moishezon. Let
\begin{align}
	\phi_\pi : \pi_1(S \setminus \Delta, b) \to \Mod(\pi^{-1}(b)) \cong \Mod(\Sigma_1) \cong \SL_2\ZZ \label{defn:mon-rep}
\end{align}
be the monodromy representation of $\pi$, choosing some identification $\pi^{-1}(b)$ with the standard torus $\Sigma_1$ once and for all. As discussed in the introduction, the action of $\Mod(S^2,\Delta)$ acts on the collection of conjugacy classes of all such representations via its outer action on $\pi_1$. The induced action is the familiar \textit{Hurwitz action} from the theory of Lefschetz fibrations, and can be phrased in terms of the monodromy factorization. The orbit of $[\phi_\pi]$ under the Hurwitz action is referred to as the \textit{Hurwitz orbit}. Equipped with this notation, we may state the theorem of Moishezon.
\begin{theorem}[Moishezon {\cite[Lemma 7a]{moishezon}}]\label{thm:moishezon-matsumuto}
	Let $\pi : M \to S^2$ be an elliptic fibration, then $\Br(\pi)$ is the stabilizer of $[\phi_\pi] \in \mathfrak{X}_\ZZ(S^2,\Delta)$, i.e, the conjugacy class of the monodromy representation $\phi_\pi$, under the action of $\Mod(S^2,\Delta)$. 
\end{theorem}

\begin{remark}
	Moishezon original phrasing of \Cref{thm:moishezon-matsumuto} is slightly different. However, the statement in \Cref{thm:moishezon-matsumuto} follows immediately. See Endo's survey for a similar statement in the literature \cite[Theorem 3.3]{endo}.
\end{remark}

Using \Cref{thm:moishezon-matsumuto}, Moishezon classified elliptic fibrations by showing that all conjugacy classes of homomorphisms $\pi_1(S^2 \setminus \Delta) \to \Mod(\Sigma_1) \cong \SL_2\ZZ$ that take simple loops about the punctures to Dehn twists lie in a single Hurwitz orbit. Note that the presence of a single such homomorphism implies that $n = \abs{\Delta}$ is a multiple of $12$, since $(\SL_2\ZZ)^{ab} = \ZZ/12\ZZ$ and each Dehn twist is sent under the abelianization map to $1 \in \ZZ/12\ZZ$.
\begin{theorem}[Moishezon {\cite[Theorem 9]{moishezon}}]\label{thm:moishezon}
	The number $n = 12d$ of singular fibers of a genus one Lefschetz fibration $\pi : M \to S^2$ determines $\pi$ in the following sense: given two genus one fibrations $\pi, \pi' : M,M' \to S^2$ with the same number of singular fibers, there are diffeomorphisms $F,f$ making the following diagram commute
	\begin{center}
		\begin{tikzcd}
			M \ar[d] \ar[r,"F"] & M' \ar[d] \\
			(S^2,\Delta_\pi) \ar[r,"f"] & (S^2, \Delta_{\pi'}).
		\end{tikzcd}
	\end{center}
	Furthermore, simple loops $\gamma_1,\ldots,\gamma_n$ on $(S^2,\Delta)$ may be chosen so that the monodromy factorization of $\pi$ is
	\begin{align*}
		(T_\alpha T_\beta)^{6d} = \Id,
	\end{align*}
	where $\alpha,\beta$ are curves with geometric intersection number $i(\alpha,\beta) = 1$ and $T_\alpha,T_\beta$ denote the respective Dehn twists about these curves.
\end{theorem}
In general, following Endo, we call two genus one fibrations $\pi,\pi' : M,M' \to S^2$ \textit{weakly isomorphic} provided there exists such a pair $(F,f)$ of diffeomorphisms so that
\begin{center}
	\begin{tikzcd}
		M \ar[d] \ar[r,"F"] & M' \ar[d] \\
		(S^2,\Delta_\pi) \ar[r,"f"] & (S^2, \Delta_{\pi'})
	\end{tikzcd}
\end{center}
commutes \cite[Definition 3.1]{endo}. Changing $\pi$ by a weak isomorphism changes $\Br(\pi)$ by conjugation in $\Mod(S^2,\Delta_\pi)$. Furthermore $\Mod(\pi)$, as defined in the introduction, consists of the isotopy classes of weak isomorphisms from $\pi$ to itself.

\subfile{clean-proof.tex}

\subfile{orders-not-n.tex}

\section{Difficulties Inherent in Classifying Small Orders}\label{sec:distinct-orders}
We now explain why the methods above do not generalize to classifying elements of order $n/a$ for $a \mid n$. As an example of an order $n/2$ element of $\Br(\pi)$, one may consider simple closed curves $\alpha,\beta,\gamma,\de$ forming the vertices of two triangles $(\alpha,\gamma,\beta)$ and $(\alpha,\beta,\de)$ labeled clockwise in the Farey complex (see \Cref{fig:farey-fin-order}). The monodromy $\phi$ corresponding to the factorization
\begin{align}
	((T_\alpha T_\gamma T_\beta)(T_\beta T_\de T_\alpha))^{n/6} = \Id \label{eq:order-n2}
\end{align}
is not conjugate to $r \cdot \phi$, where $r$ is the order $n$ rotation. However it is conjugate to $r^3 \cdot \phi$. Hence the monodromy factorization \eqref{eq:order-n2} corresponds to a ``new'' order $n/3$ conjugacy class not arising as the cube of either of the order $n$ conjugacy classes discussed in \Cref{thm:fin-order}, and so the order $n/2$ element corresponding to \eqref{eq:order-n2} does not belong to a $\ZZ/n\ZZ$ subgroup of $\Br(\pi)$.

\begin{figure}[h]
	\centering
	\includestandalone[width=.35\textwidth]{tikz/farey-complex}
	\caption{Farey Tesselation labeled $\alpha,\de,\beta$ as well as points $i,\omega = e^{2\pi i/3}$ under the identification of the Poincar\'{e} disk with $\mathbb{H}^2$.}
	\label{fig:farey-fin-order}
\end{figure}

\begin{figure}[ht]
	\centering
	\subfloat[$w_2 = 0$]{
		\centering
		\includegraphics[scale=0.17]{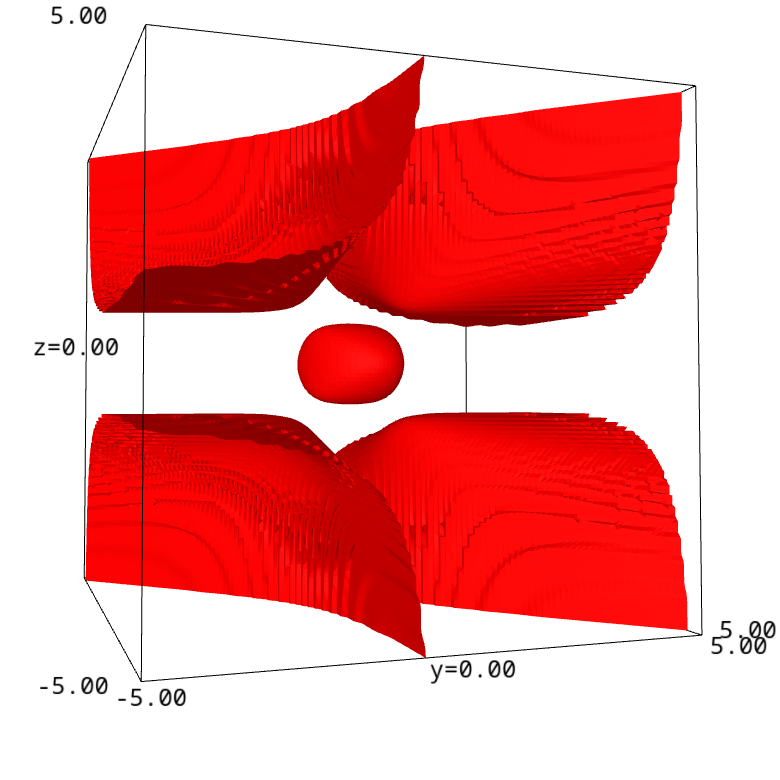}
	} 
	\hspace{0.1in} 
	\subfloat[$w_2 = -4$]{
		\centering
		\includegraphics[scale=0.17]{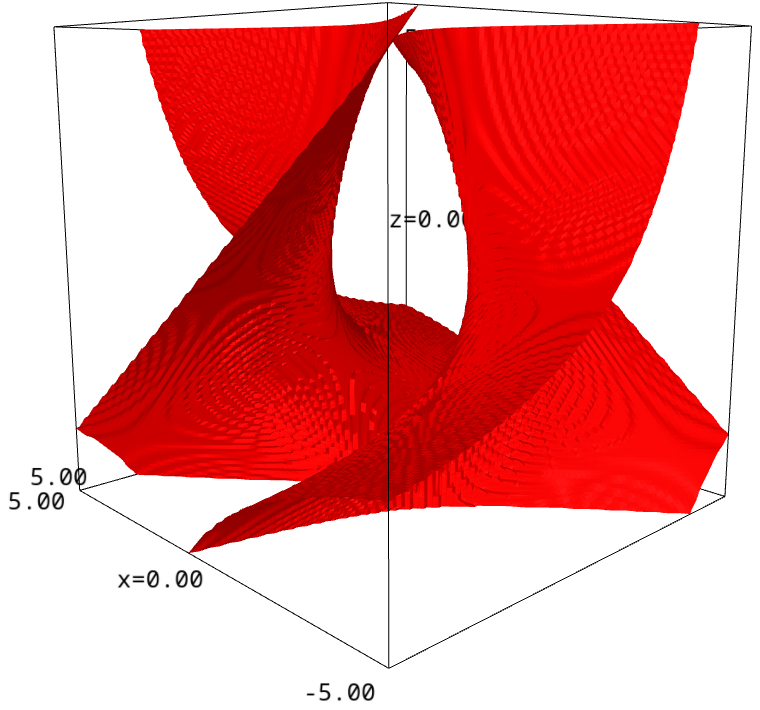}
	}\\
	\caption{Noncompact 3-dimensional slices of $g_6 = 0$ where $w = (w_1, w_2)$ is restricted to a hyperplane.}
	\label{fig:noncompact-fin-order}
\end{figure}

When $a = 2$, a similar argument to the proofs of \Cref{thm:fin-order,thm:matrix-fin-order} shows that classifying elements of order $n/2$ is essentially equivalent to classifying primitive vectors $v,w \in \ZZ^2$ so that
\begin{align*}
	\abs{\tr T_v T_w T_{Cv} T_{Cw} \cdots T_{C^{n/2-1}v}T_{C^{n/2-1}w}} \leq 2.
\end{align*}
As before, $C$ is of finite order, and so can be conjugated to one of $C_3,C_4,C_6$. Thus we obtain three polynomials $g_3(v,w), g_4(v,w), g_6(v,w)$, and we wish to find the lattice points in $\ZZ^4$ satisfying $\abs{g_j(v,w)} \leq 2$. Two three-dimensional slices of $g_6 = 0$ are displayed in \Cref{fig:noncompact-fin-order}. 

There are two primary obstructions to generalizing the proof of \Cref{thm:fin-order} above, which are implicitly linked:
\begin{enumerate}
	\item The invariant ring $\ZZ[v,w]^H$, where $H$ is the centralizer of $C$ acting diagonally on the $v,w$ coordinates, is significantly more complex than the invariant ring computed in \Cref{prop:invariant-theory}.
	\item The real algebraic locus of $\abs{g_j(v,w)} \leq 2$ in $\R^4$ is not compact, as can be seen in \Cref{fig:noncompact-fin-order}.
\end{enumerate}
Similar difficulties arise when considering which powers of $\tau$ and $\eta$ lie in some conjugate of $\Br(\pi)$, where $\tau$ and $\eta$ are the order $n-1$ and $n-2$ rotation respectively. For example, consider the following monodromy factorization for an elliptic fibration $\pi$ with $24$ singular fibers, i.e., of a K3 surface:
\begin{align}
	T_\alpha (T_\beta T_\alpha)^6 (T_\beta T_\alpha)^5 T_\beta = \Id,\label{eq:eta-12}
\end{align}
where $\alpha,\beta$ are curves with intersection number $i(\alpha,\beta) = 1$. The monodromy factorization \eqref{eq:eta-12} is fixed by $\eta^{12}$, which can be verified using the generating set for $\pi_1(S^2\setminus \Delta)$ given in \Cref{sec:orders-not-n}. Hence \eqref{eq:eta-12} represents an element of order $11$ in $\Br(\pi)$, while there is no element of order $22$, contrasting sharply with Murasugi's theorem. Despite these difficulties, the symmetry involved in the construction of the monodromy factorizations corresponding to the two order $n$ conjugacy classes and the conjugacy class of order $n/2$ displayed in \eqref{eq:order-n2} suggests leveraging the geometry of the Farey complex as a general approach.

\printbibliography

\end{document}

%% file: research-notes-preamble.tex
\usepackage{subfiles}

\usepackage{amsfonts,amsthm,amssymb,amsmath,amscd,euscript}
\usepackage{nicematrix}
\usepackage{mathdots}
\usepackage{csquotes}
\usepackage{mathrsfs}
\usepackage{wasysym}
\usepackage{bbm}
\usepackage{ifthen}
\usepackage{xpatch}
\usepackage{thmtools}
\usepackage{thm-restate}
\usepackage{calc}
\usepackage[reftex]{theoremref}

\usepackage[backend=biber,style=alphabetic,giveninits=true]{biblatex}
\usepackage[english]{babel}
\usepackage[hypertexnames=false]{hyperref}

\usepackage{multicol}
\usepackage[dvipsnames,rgb]{xcolor}
\usepackage{framed}
\usepackage{pdfpages}
\usepackage{scalerel}
\usepackage{fullpage}
\usepackage{color}
\usepackage{lipsum}
\usepackage{bold-extra}
\usepackage{fancyhdr}

\usepackage{tikz,tikz-3dplot}
\usepackage{tkz-berge}

\usepackage{pgfplots}
\pgfplotsset{compat=1.7}
\usetikzlibrary{shapes, arrows}
\usetikzlibrary{shapes.geometric}
\usetikzlibrary{graphs}
\usetikzlibrary{graphs.standard}
\usetikzlibrary{decorations.pathmorphing}
\usetikzlibrary{decorations.markings}
\tikzset{->-/.style={decoration={
  markings,
  mark=at position #1 with {\arrow{>}}},postaction={decorate}},
  ->-/.default=0.53
}
\tikzset{-<-/.style={decoration={
  markings,
  mark=at position #1 with {\arrow{<}}},postaction={decorate}},
  -<-/.default=0.53
}
\pgfplotsset{soldot/.style={color=blue,only marks,mark=*}} \pgfplotsset{holdot/.style={color=blue,fill=white,only marks,mark=*}}
\usepackage{tikz-cd}
\usepackage[all,2cell,cmtip]{xy}
\UseTwocells

\usepackage{mathtools}
\usepackage{float}
\usepackage[margin=1in,heightrounded]{geometry}
\usepackage{etoolbox}
\usepackage[most]{tcolorbox}
\tcbuselibrary{theorems}
\tcbuselibrary{skins}
\usepackage{mdframed}
\usepackage[inline, shortlabels]{enumitem}
\SetEnumitemKey{twocol}{itemsep = 1\itemsep,parsep  = 1\parsep,before  = \raggedcolumns\begin{multicols}{2},after   = \end{multicols}}
\usepackage{tasks}
\usepackage{stackengine}
\usepackage{tabularx}
\usepackage{soul}
\usepackage[a]{esvect}
\usepackage[activate={true,nocompatibility},nopatch=footnote,final,tracking=true,kerning=true,spacing=true,factor=1100,stretch=10,shrink=10]{microtype}
\microtypecontext{spacing=nonfrench}
\usepackage{siunitx}
\usepackage{mathrsfs}
\usepackage{quiver} % https://q.uiver.app/quiver.sty
\usepackage{diagbox}
\DeclareFontFamily{U}{mathx}{\hyphenchar\font45}
\DeclareFontShape{U}{mathx}{m}{n}{
      <5> <6> <7> <8> <9> <10>
      <10.95> <12> <14.4> <17.28> <20.74> <24.88>
      mathx10
      }{}
\DeclareSymbolFont{mathx}{U}{mathx}{m}{n}
\DeclareFontSubstitution{U}{mathx}{m}{n}
\DeclareMathAccent{\widecheck}{0}{mathx}{"71}
\DeclareMathAccent{\wideparen}{0}{mathx}{"75}

\usepackage[color=red!60, colorinlistoftodos]{todonotes}
\usepackage{setspace}
\usepackage{stmaryrd}
\usepackage{centernot}
\usepackage{subfig}
\usepackage{cancel}
\newcommand\Ccancel[2][black]{
    \let\OldcancelColor\CancelColor
    \renewcommand\CancelColor{\color{#1}}
    \cancel{#2}
    \renewcommand\CancelColor{\OldcancelColor}
}
\usepackage{makecell}
\hypersetup{colorlinks=true,citecolor=blue,urlcolor =blue,linkbordercolor={1 0 0}}

\definecolor{bluepigment}{rgb}{0.2, 0.2, 0.6}
\definecolor{cobalt}{rgb}{0.0, 0.28, 0.67}
\definecolor{darkpowderblue}{rgb}{0.0, 0.2, 0.6}
\definecolor{egyptianblue}{rgb}{0.06, 0.2, 0.65}
\definecolor{lapislazuli}{rgb}{0.15, 0.38, 0.61}
\definecolor{tealblue}{rgb}{0.21, 0.46, 0.53}
\definecolor{darktealblue}{rgb}{0.11, 0.36, 0.43}
\definecolor{zaffre}{rgb}{0.0, 0.08, 0.66}
\definecolor{zaffreprime}{rgb}{0.4, 0.18, 0.76}
\definecolor{arsenic}{rgb}{0.23, 0.27, 0.29}
\definecolor{bulgarianrose}{rgb}{0.28, 0.02, 0.03}
\definecolor{aurometalsaurus}{rgb}{0.43, 0.5, 0.5}
\definecolor{amethyst}{rgb}{0.6, 0.4, 0.8}

\definecolor{pblue}{rgb}{.114, .102, .365}
\definecolor{pblue2}{rgb}{.204, .196, .471}
\definecolor{pblue3}{rgb}{.655,.784,.855}
\definecolor{pgreen}{rgb}{.075,.212,.337}
\definecolor{pgreen3}{rgb}{.667,.882,.765}

\definecolor{codegreen}{rgb}{0,0.6,0}
\definecolor{codegray}{rgb}{0.5,0.5,0.5}
\definecolor{codepurple}{rgb}{0.58,0,0.82}
\definecolor{backcolour}{rgb}{0.95,0.95,0.92}

\allowdisplaybreaks[1]

\makeatletter
\patchcmd{\@maketitle}%
  {\ifx\@empty\@dedicatory}%
  {\ifx\@empty\@date \else {\vskip3ex \centering\footnotesize\@date\par\vskip1ex}\fi%
   \ifx\@empty\@dedicatory}%
  {}{}%
\patchcmd{\@adminfootnotes}%
  {\ifx\@empty\@date\else \@footnotetext{\@setdate}\fi}%
  {}{}{}%

\def\section{\@startsection{section}{1}%
  \z@{.7\linespacing\@plus\linespacing}{.5\linespacing}%
  {\normalfont\Large\bfseries}}%{\normalfont\Large\bfseries\scshape\centering}} % full range of options tested
\def\subsection{\@startsection{subsection}{2}%
  \z@{.5\linespacing\@plus.7\linespacing}{.5\linespacing}%
  {\normalfont\large\bfseries}}%
\def\subsubsection{\@startsection{subsubsection}{3}%
  \z@{.5\linespacing\@plus.7\linespacing}{.5\linespacing}%
  {\normalfont\bfseries}}%
\def\@tocline#1#2#3#4#5#6#7{\relax%
  \ifnum #1>\c@tocdepth % then omit
  \else%
    \par \addpenalty\@secpenalty\addvspace{#2}%
    \begingroup \hyphenpenalty\@M%
    \@ifempty{#4}{%
      \@tempdima\csname r@tocindent\number#1\endcsname\relax%
    }{%
      \@tempdima#4\relax%
    }%
    \parindent\z@ \leftskip#3\relax \advance\leftskip\@tempdima\relax%
    \rightskip\@pnumwidth plus4em \parfillskip-\@pnumwidth%
    #5\leavevmode\hskip-\@tempdima%
      \ifcase #1%
       \or\or \hskip 1em \or \hskip 2em \else \hskip 3em \fi%
      #6\nobreak\relax%
    \dotfill\hbox to\@pnumwidth{\@tocpagenum{#7}}\par%
    \nobreak%
    \endgroup%
  \fi}%
\makeatother%

\makeatletter%
\xpatchcmd{\endmdframed}%
    {\aftergroup\endmdf@trivlist\color@endgroup}%
    {\endmdf@trivlist\color@endgroup\@doendpe}%
    {}{}%
\makeatother%

\newenvironment{solution} {\begin{mdframed}[skipabove=0.1in,skipbelow=0.1in]\begin{proof}[Solution]} {\end{proof}\end{mdframed}}

\mdfdefinestyle{theoremstyle}{fontcolor=amethyst, linecolor=amethyst, bottomline=false,topline=false,rightline=false,skipabove=2pt,skipbelow=2pt}

\theoremstyle{theorem}

\declaretheorem[numberwithin=section, style=theorem, name=Theorem]{theorem}

\declaretheorem[ style=theorem, name=Lemma, sibling=theorem]{lemma}

\declaretheorem[sibling=theorem, style=theorem, name=Proposition]{proposition}

\declaretheorem[,numberwithin=section, style=theorem, name=Remark,sibling=theorem,style=definition]{remark}

\declaretheorem[sibling={example}, style=theorem, name=Exercise]{exercise}

\declaretheorem[sibling={example}, style=theorem, name=Problem]{prob}
\declaretheorem[sibling={example}, style=theorem, name=Question]{quest}

\declaretheorem[numbered=no, style=theorem, name=Problem]{prob-nonum}

\lstdefinestyle{mystyle}{
    keywordstyle=\color[rgb]{0,0,1},
    commentstyle=\color[rgb]{0.026,0.412,0.595},
    stringstyle=\color[rgb]{0.627,0.126,0.941},
    numberstyle=\color[rgb]{0.205, 0.142, 0.73},
    frame=single,
    breakatwhitespace=false,
    breaklines=true,
    postbreak=\mbox{\textcolor{darkgray}{$\hookrightarrow$}\space},
    captionpos=b,
    keepspaces=true,
    numbers=left,
    numbersep=5pt,
    showspaces=false,
    showstringspaces=false,
    showtabs=false,
    tabsize=2,
    basicstyle=\ttfamily,
    columns=fullflexible,
    keepspaces=true,
    language=C++
}

\DeclareRobustCommand{\SkipTocEntry}[3]{}
\usepackage{cleveref}

\crefname{defn}{definition}{definitions}
\Crefname{defn}{Definition}{Definitions}
\crefname{theorem}{theorem}{theorems}
\Crefname{theorem}{Theorem}{Theorems}
\crefname{prob}{problem}{problems}
\Crefname{prob}{Problem}{Problems}
\crefname{exercise}{exercise}{exercises}
\Crefname{exercise}{Exercise}{Exercises}
\crefname{proposition}{proposition}{propositions}
\Crefname{proposition}{Proposition}{Propositions}
\makeatletter
\newcommand{\leqnomode}{\tagsleft@true\let\veqno\@@leqno}
\makeatother

\newcommand{\ZZ}{\mathbb{Z}}
\newcommand{\R}{\mathbb{R}}

\newcommand{\CC}{\mathbb{C}}

\newcommand{\abs}[1]{\left\lvert #1 \right\rvert}

\newcommand{\0}{\emptyset}

\renewcommand*{\st}{\;\vert\;}

\DeclareMathOperator{\Id}{Id}

\DeclareMathOperator{\im}{im}

\DeclareMathOperator{\Hom}{Hom}

\newcommand{\rest}[2]{{% we make the whole thing an ordinary symbol
  \left.\kern-\nulldelimiterspace % automatically resize the bar with \right
  #1 % the function
  \vphantom{\big|} % pretend it's a little taller at normal size
  \right|_{#2} % this is the delimiter
  }}

\DeclareMathOperator{\GL}{GL}
\DeclareMathOperator{\SL}{SL}

\DeclareMathOperator{\tr}{tr}

\newcommand{\de}{\delta}
\renewcommand*{\d}{\mathop{}\!\mathrm{d}}

\DeclareMathOperator{\Mod}{Mod}

\DeclareMathOperator{\Stab}{Stab}

\DeclareMathOperator{\Diff}{Diff}

\newcommand*{\PP}{\mathbb{P}}

%% file: intro.tex
A \emph{smooth elliptic fibration} $\pi : M \to S^2$ of an oriented $4$-manifold $M$ is a smooth map $M \to S^2$ with finitely many critical values $\Delta \subseteq S^2$ such that the fibers $\pi^{-1}(b)$ are smooth elliptic curves for $b \in S^2 \setminus \Delta$ and nodal elliptic curves for $b \in \Delta$. Note that this equips $\pi$ with a section{\footnotemark} $s : S^2 \to M$ picking out the basepoint of each elliptic curve. 
{\footnotetext{For us, we consider only elliptic fibrations $\pi : M \to S^2$ with nodal singularities without multiples. Further, the role of a section is auxiliary to our results, and so one may instead consider genus one fibrations.}}
Examples include the rational elliptic surface, given by blowing up $\PP^2$ along the intersection points of two generic cubics, and the elliptically fibered K3 surfaces. Given an elliptic fibration
$\pi : M \to S^2$, let
\begin{align*}
   \Diff^+(\pi) &\coloneqq \{F \in \Diff^+(M) \st F \text{ takes fibers of } \pi \text{ to fibers of } \pi\}.
\end{align*}
The \textit{smooth automorphism group} of $\pi$ as defined by Farb--Looijenga (see \cite{fl-smooth-mw}) is 
\begin{align*}
   \Mod(\pi) &\coloneqq \pi_0(\Diff^+(\pi)),
\end{align*}
Tracking the location of the singular fibers under an element of $\Diff^+(\pi)$ gives a \emph{braid monodromy representation}:
\begin{align*}
   \rho : \Mod(\pi) \to \Mod(S^2, \Delta),
\end{align*}
where $\Delta \subseteq S^2$ is the singular locus of $\pi$ and $\Mod(S^2,\Delta) \coloneqq \pi_0(\Diff^+(S^2,\Delta))$, with $\Diff^+(S^2,\Delta)$ consisting of diffeomorphisms of $S^2$ mapping $\Delta$ to itself setwise.  We study the image
\begin{align}
   \Br(\pi) \coloneqq \im(\rho) < \Mod(S^2, \Delta), \label{defn:br-pi}
\end{align}
of the braid monodromy $\rho$. Thus $\Br(\pi)$ associates to each elliptic fibration $\pi$ a subgroup of the spherical braid group{\footnotemark} $\Mod(S^2,\Delta)$. 
\footnotetext{There are two distinct notions of the spherical braid group, one is $B_n(S^2) \coloneqq \pi_1(\operatorname{Conf}_n(S^2))$, the fundamental group of $n$-point configurations in $S^2$, the other is $\Mod(S^2, \{p_1,\ldots,p_n\})$, a marked mapping class group. These differ by an exact sequence $1 \to \ZZ/2\ZZ \to B_n(S^2) \to \Mod(S^2, \{p_1,\ldots,p_n\}) \to 1$ \cite[\S 9.1.4]{primer}. In this paper we restrict ourselves to the latter notion.}
The subgroup $\Br(\pi) < \Mod(S^2, \Delta)$ consists of those braids that lift to a fiber-preserving diffeomorphism of $M$. For convenience, let $n \coloneqq \abs{\Delta}$ be the number of singular fibers. In this paper, we classify the conjugacy classes of elements of order $n$ in $\Br(\pi)$ and show that there are no elements of order $n-1$ or $n-2$, as there are in $\Mod(S^2,\Delta)$.

Classifying torsion elements in $\Br(\pi)$ up to $\Br(\pi)$ conjugacy is intimately tied to a theorem of Murasugi which classifies the finite order elements of $\Mod(S^2,\Delta)$ up to conjugacy in $\Mod(S^2,\Delta)$ \cite{murasugi-braid-torsion}. Murasugi shows that every finite order element of $\Mod(S^2,\Delta)$ is conjugate to a power of one of the following:
\begin{enumerate}
   \item The order $n$ rotation $\sigma_1 \cdots \sigma_{n-1}$.
   \item The order $n-1$ rotation, achieved by placing one marked point at the north pole and the remainder along the equator.
   \item The order $n-2$ rotation, achieved by placing one marked point at the north pole, another at the south pole, and the remainder along the equator.
\end{enumerate}
A similar classification in $\Br(\pi)$ does not follow formally, however. First, each conjugate $bgb^{-1}$ of a finite order element $g \in \Mod(S^2,\Delta)$ for $b \in \Mod(S^2,\Delta)$ may or may not appear in the subgroup $\Br(\pi) < \Mod(S^2,\Delta_\pi)$. The conditions for whether such an element $bgb^{-1}$ lies in $\Br(\pi)$ are given in \Cref{sec:prelim} and further explained in \Cref{sec:equiv}. Second, a single conjugacy class in $\Mod(S^2,\Delta_\pi)$ may split into multiple distinct conjugacy classes in $\Br(\pi)$. The relationship between these classifications is further complicated by the fact that the index $[\Mod(S^2,\Delta_\pi) : \Br(\pi)] = \infty$, which we show in \cite{braid-index-j}. Nonetheless, we prove the following.
\begin{theorem}[\textbf{Torsion in $\Br(\pi)$ of order $n, n-1,$ and $n-2$}]\label{thm:fin-order} Let $\pi : M \to S^2$ be an elliptic fibration with $n$ nodal fibers. The following hold:
	\begin{enumerate}[(a)]
		\item\label{item:ordern} There are exactly two distinct conjugacy classes in $\Br(\pi)$ of elements with order $n$. Furthermore, up to replacing $\Br(\pi)$ by a conjugate subgroup in $\Mod(S^2,\Delta)$, these two classes are represented by $r = \sigma_1\cdots \sigma_{n-1} \in \Br(\pi)$ and $r^{-1} = T_nrT_n^{-1} \in \Br(\pi)$, where $T_n$ is the \textit{Garside half-twist} defined by 
			 \begin{align*}
  			 T_n	= (\sigma_1\cdots \sigma_{n-1})(\sigma_1\cdots\sigma_{n-2})\cdots (\sigma_1\sigma_2)\sigma_1,
  			 \end{align*}
  			 and $\sigma_1,\ldots,\sigma_{n-1}$ are standard half-twist generators of $\Mod(S^2,\Delta_\pi)$.
		 \item\label{item:order-n-minus-one} There are no elements of order $n-1$ in $\Br(\pi)$,
		 \item\label{item:order-n-minus-two} There are no elements of order $n-2$ in $\Br(\pi)$.
	\end{enumerate}
\end{theorem}
Note that, despite part (c) of \Cref{thm:fin-order} ruling out elements of order $n-2$, there do exist elements whose orders divide $n-2$ and do not divide $n-1$ or $n$ itself. For an example, see \eqref{eq:eta-12} in \Cref{sec:distinct-orders}. The proof of \Cref{thm:fin-order} relies on an understanding of the Hurwitz action of $\Mod(S^2,\Delta)$ on the $\SL_2$-character variety for $(S^2,\Delta)$, which we now recall.

\medskip
\para{$\bm{\Br(\pi)}$ and $\bm{\SL_2}$-character varieties}

We prove \Cref{thm:fin-order} by showing that each part is equivalent to a corresponding statement classifying simultaneous conjugacy classes of $n$-tuples of matrices in $\SL_2\ZZ$ which satisfy particular conditions. 
We first detail the relationship between $\Br(\pi)$ and the Hurwitz action on the $\SL_2$-character variety.

Let $\pi : M \to S^2$ be an elliptic fibration. A choice of basepoint $b \in S^2 \setminus \Delta$ gives an associated monodromy representation
\[
	\phi_{\pi} : \pi_1(S^2 \setminus \Delta,b) \to \Mod(\pi^{-1}(b)) \cong \Mod(\Sigma_1) \cong \SL_2\ZZ,
\]
by identifying $\pi^{-1}(b)$ with the standard torus $\Sigma_1$. Note that changing the basepoint $b$ or the identification of $\pi^{-1}(b)$ with $\Sigma_1$ changes $\phi_\pi$ by conjugation, and hence associated to $\pi$ is a point $[\phi_\pi]$ in the \textit{integral character variety}\footnotemark
\footnotetext{$\mathfrak{X}_\ZZ(S^2,\Delta)$ itself is not a variety. However, the GIT quotient $\Hom(\pi_1(S^2 \setminus \Delta),\SL_2\CC)\sslash \SL_2\CC$) does form a variety, commonly referred to as the character variety. We refer to $\mathfrak{X}_\ZZ(S^2,\Delta)$ as the integral character variety for convenience and for context.}
\begin{align*}
	\mathfrak{X}_\ZZ(S^2,\Delta) \coloneqq \Hom(\pi_1(S^2 \setminus \Delta), \SL_2\ZZ)/\SL_2\ZZ,
\end{align*}
where $\SL_2\ZZ$ acts by simultaneous conjugation. The spherical braid group $\Mod(S^2,\Delta)$ naturally acts on $\mathfrak{X}_\ZZ(S^2,\Delta)$ via precomposition by (outer) automorphisms on $\pi_1$. The induced action of $\Mod(S^2,\Delta)$ on $\mathfrak{X}_\ZZ(S^2,\Delta)$ is referred to as the \textit{Hurwitz action}. A theorem of Moishezon (see \Cref{thm:moishezon-matsumuto} below) implies that the liftable braids are precisely the stabilizers of $[\phi_\pi] \in \mathfrak{X}_\ZZ(S^2,\Delta)$; that is
\begin{align}
	\Br(\pi) = \Stab_{\Mod(S^2,\Delta)} [\phi_\pi].\label{eq:br-pi-stab}
\end{align}
We discuss this further in \Cref{sec:prelim}. For now, recall that choosing generators $\gamma_1,\ldots,\gamma_n$ for $\pi_1(S^2,\Delta)$ surrounding each puncture counterclockwise realizes the representation $\phi_\pi$ as a \textit{factorization} in the mapping class group
\begin{align*}
	\phi_\pi(\gamma_1)\cdots \phi_\pi(\gamma_n) = \Id.
\end{align*}
By Picard--Lefschetz theory, each $\phi_\pi(\gamma_i)$ is a Dehn twist about a simple closed curve in $\pi^{-1}(b) \cong \Sigma_1$ \cite{endo}. We refer to the ordered tuple $(\phi_\pi(\gamma_i)) \in (\SL_2\ZZ)^n$ as the \textit{monodromy factorization} associated to the fibration $\pi$. 

We now return to \Cref{thm:fin-order}. By Murasugi's theorem, an order $n$ element of $\Br(\pi)$ is given as a conjugate $brb^{-1} \in \Br(\pi)$ with $b \in \Mod(S^2,\Delta)$, where $r$ is the order $n$ rotation. Applying \eqref{eq:br-pi-stab} shows that $rb^{-1} \cdot [\phi_\pi] = b^{-1} \cdot [\phi_\pi]$, i.e., that these representations are conjugate.
We call $b^{-1} \cdot [\phi_\pi]$ a \textit{rotation invariant monodromy representation}.
By choosing generators $\gamma_i \in \pi_1(S^2 \setminus \Delta)$ appropriately, we obtain a notion of a \textit{rotation invariant monodromy factorization}.
In \Cref{sec:equiv}, we show that part (a) of \Cref{thm:fin-order} is equivalent to the following theorem classifying rotation invariant monodromy factorizations up to conjugacy.
\begin{theorem}[\textbf{Rotation Invariant Monodromy Factorizations}]\label{thm:matrix-fin-order}
Call an ordered $n$-tuple $(X_1,\ldots,X_n)$ in $(\SL_2\ZZ)^n$ a \textit{rotation invariant monodromy factorization} provided that for some $C \in \SL_2\ZZ$:
	\begin{enumerate}
		\item\label{item:prim} $X_1 = T_v$, where $T_v(x) = x + (x,v)v$ is a symplectic transvection in a primitive vector $v \in \ZZ^2$. Equivalently, $X_1$ is conjugate to $\begin{psmallmatrix} 1 & -1 \\ 0 & 1 \end{psmallmatrix} $.
		\item\label{item:cyclic} $X_{i+1} = CX_iC^{-1}$ for $i \in \ZZ/n\ZZ$
		\item\label{item:prod-is-id} $\prod_{i=1}^n X_i = \Id$
	\end{enumerate}
	Then there are exactly two rotation invariant monodromy factorizations $(X_1,\ldots,X_n)$ in $(\SL_2\ZZ)^n$ up to simultaneous conjugation.
	Explicitly, for any such sequence $(X_1,\ldots,X_n)$, $n$ is divisible by $12$ and there is some matrix $D \in \SL_2\ZZ$ so that exactly one of the following holds
	\begin{enumerate}[(a)]
		\item \quad\quad $
				(DX_1D^{-1},\ldots,DX_n D^{-1}) = \left(\begin{pmatrix} 1 & -1 \\ 0 & 1 \end{pmatrix}, \begin{pmatrix} 1 & 0 \\ 1 & 1 \end{pmatrix}, \begin{pmatrix} 1 & -1 \\ 0 & 1 \end{pmatrix}, \ldots \right)
			$ or
		\item \quad\quad $
				(DX_1D^{-1},\ldots,DX_n D^{-1}) = \left(\begin{pmatrix} 1 & -1 \\ 0 & 1 \end{pmatrix}, \begin{pmatrix} 2 & -1 \\ 1 & 0 \end{pmatrix}, \begin{pmatrix} 1 & 0 \\ 1 & 1 \end{pmatrix}, \begin{pmatrix} 1 & -1 \\ 0 & 1 \end{pmatrix}, \ldots \right)$
	\end{enumerate}
	The tuples falling into case (a) are said to have period 2 and the tuples falling into case (b) are said to have period $3$.
\end{theorem}
Experts may notice the similarity between \Cref{thm:matrix-fin-order} and a theorem of Moishezon and Seiler, which states that any two tuples $(X_1,\ldots, X_n)$ satisfying $\prod_{i=1}^n X_i = \Id$ and with $X_i \in \SL_2\ZZ$ conjugate to $\begin{pmatrix} 1 & -1 \\ 0 &1 \end{pmatrix}$ lie in the same Hurwitz orbit (up to simultaneous conjugation) \cite[Lemma 8]{moishezon}. However, in contrast, \Cref{thm:matrix-fin-order} states that there are only two rotation invariant monodromy factorizations up to simultaneous conjugacy within the Hurwitz orbit.

We delay the precise statements of the analogous classifications corresponding to parts (b) and (c) of \Cref{thm:fin-order} to \Cref{sec:orders-not-n}, given there as as \Cref{thm:matrix-fin-order-n-minus-one} and \Cref{prop:matrix-fin-order-n-minus-two}.

\medskip
\para{Geometric interpretation of $\bm{r}$ and $\bm{r^{-1} = T_n rT_n^{-1}}$}

Let $\pi : M \to S^2$ be an elliptic fibration with $n$ singular fibers. Observing that $(\SL_2\ZZ)^{ab} = \ZZ/12\ZZ$ and that all Dehn twists are conjugate in $\SL_2\ZZ$, we find that $n = 12d$ for some $d \geq 1$. The two conjugacy classes $r$ and $r^{-1} = T_nrT_n^{-1}$ in $\Br(\pi)$ are associated to two distinct monodromy factorizations for $\pi$, given by different choices of simple closed curves $\gamma_i$ (resp. $\gamma_i'$) enclosing each puncture. In \Cref{thm:matrix-fin-order}, $\rho_\pi(\gamma_i) \in \SL_2\ZZ$ (resp. $\rho_\pi(\gamma_i')$) appear as the matrices given in (a) and (b). Interpreted geometrically, these factorizations are given by
\begin{align}
	(T_\alpha \cdot T_\beta)^{6d} = \Id && (T_\alpha \cdot T_\de \cdot T_\beta)^{4d} = \Id, \label{eq:fin-order-mon}
\end{align}
respectively, where the curves $\alpha,\beta,\de$ are the vertices of a triangle labeled counterclockwise on the Farey tesselation (see \Cref{fig:farey-fin-order} in \Cref{sec:distinct-orders} below) and $T_\alpha,T_\beta,T_\de$ are the corresponding Dehn twists. Each of these monodromy factorizations are fixed up to conjugacy by $r$ under the Hurwitz action, and the latter is obtained from $(T_\alpha \cdot T_\beta)^{6d}$ by applying the Garside half-twist $T_n$. We can associate the first monodromy factorization to the order 4 automorphism of the elliptic curve $E_i = \CC/\ZZ\langle 1,i \rangle$, and the second monodromy factorization to the order 6 automorphism of the elliptic curve $E_\omega = \CC/\ZZ\langle 1,\omega\rangle$ with $\omega = e^{2\pi i/3}$. The automorphisms of these elliptic curves naturally appear in the proof of \Cref{thm:fin-order}. Unsurprisingly, the number rings $\ZZ[\omega]$ and $\ZZ[i]$ also appear and play a pivotal role.

\medskip
\para{Method of Proof for \Cref{thm:matrix-fin-order}}

To prove \Cref{thm:matrix-fin-order} we first prove that the conjugating matrix $C$ has finite order. Thus $T_vT_{Cv}\cdots T_{C^{k-1}v}$ has finite order, where $C^k = -\Id$. The second step of the proof consists of analyzing and computing the trace polynomial
\begin{align*}
	f(v) = \tr(T_vT_{Cv}\cdots T_{C^{k-1} v}).
\end{align*}
We then show that the level sets of this polynomial are ellipses, at which point we can simply identify the integral points $v \in \ZZ^2$ where $\abs{f(v)} \leq 2$. The level sets of the trace polynomials naturally correspond to integral points in $\ZZ[i]$ and $\ZZ[\omega]$ with specified norms. We display these here as \Cref{fig:ordernconics-intro}, and delay the explanation for the proof. The method of proof for \Cref{thm:matrix-fin-order-n-minus-one} and \Cref{prop:matrix-fin-order-n-minus-two} is similar and relies on analyzing the conjugating matrix $C$ as well.

\begin{figure}[ht]
	\centering
	\input{tikz/ellipses}
	\caption{The three conics $N(p+q\omega) = 1, N(p+qi) = 2$, and $N(p+q\omega) = 3$, where $N$ is the norm of the corresponding number ring.}
	\label{fig:ordernconics-intro}
\end{figure}

\medskip
\para{Organization of the Paper}

In \Cref{sec:prelim} we explain Moishezon's theorem classifying elliptic fibrations, and review a maps version of the theorem which allows us to determine when a braid lifts to $\Diff^+(\pi)$ via the monodromy representation. Using this description, we prove the equivalence of \Cref{thm:fin-order} part (a) and \Cref{thm:matrix-fin-order} in \Cref{sec:equiv}. After doing so, we directly compute the aforementioned trace polynomial and prove \Cref{thm:matrix-fin-order} in \Cref{sec:trace-poly}. Following this, we give an algebraic motivation for the appearance of the norms on $\ZZ[i]$ and $\ZZ[\omega]$ in \Cref{sec:alg-trace-poly}. In \Cref{sec:orders-not-n}, we state and prove \Cref{thm:matrix-fin-order-n-minus-one} and \Cref{prop:matrix-fin-order-n-minus-two}, which are equivalent to part (b) and part (c) of \Cref{thm:fin-order} respectively. Finally, we indicate in \Cref{sec:distinct-orders} the difficulties that orders smaller than $n,n-1,$ and $n-2$ present for our proof.

%% file: tikz/ellipses.tex
\begin{tikzpicture}[scale=1.3]
	\begin{axis}[
		xmin=-3,xmax=3,
		ymin=-3,ymax=3,
		axis equal,
		grid=both,
		anchor=center,
		major grid style={line width=.2pt,draw=gray!100},
		axis lines=middle,
		enlargelimits={abs=0.5},
		axis line style={latex-latex},
		ticklabel style={font=\tiny,fill=white},
		xtick = {-3,-2,-1,1,2,3,4},
		ytick = {-3,-2,-1,1,2,3,4},
	]
	\end{axis}
	\draw[dashed,color=blue] (0,0) circle (1.414);
	\draw[rotate around={-45:(0,0)}, color=red] (0,0) ellipse (1.414 and 2.449);
	\draw[rotate around={-45:(0,0)}, color=purple] (0,0) ellipse (0.8165 and 1.414);
	\node at (1,0) [draw, fill=black, circle, inner sep=1pt] {};
	\node at (1,1) [draw, fill=red, circle, inner sep=1pt] {};
	\node at (0,1) [draw, fill=black, circle, inner sep=1pt] {};
	\node at (2,1) [draw, fill=black, circle, inner sep=1pt] {};
	\node at (1,2) [draw, fill=black, circle, inner sep=1pt] {};
	\node at (-1,1) [draw, fill=black, circle, inner sep=1pt] {};

	\node at (-1,0) [draw, fill=black, circle, inner sep=1pt] {};
	\node at (-1,-1) [draw, fill=red, circle, inner sep=1pt] {};
	\node at (0,-1) [draw, fill=black, circle, inner sep=1pt] {};
	\node at (-2,-1) [draw, fill=black, circle, inner sep=1pt] {};
	\node at (-1,-2) [draw, fill=black, circle, inner sep=1pt] {};
	\node at (1,-1) [draw, fill=black, circle, inner sep=1pt] {};
\end{tikzpicture}

%% file: clean-proof.tex
\section{Equivalence of \texorpdfstring{\Cref{thm:fin-order}}{Theorem 1.1} part (a) and \texorpdfstring{\Cref{thm:matrix-fin-order}}{Theorem 1.2}}\label{sec:equiv}

Let $\pi : M \to S^2$ be an elliptic fibration with $n$ nodal singular fibers located along $\Delta \subseteq S^2$. \Cref{thm:moishezon-matsumuto} identifies the subgroup $\Br(\pi)$ of liftable braids with the stabilizer of the conjugacy class of the monodromy representation. This identification is the essential ingredient in showing part (a) of \Cref{thm:fin-order} and \Cref{thm:matrix-fin-order} are equivalent. Let
\[
	\phi_\pi : \pi_1(S^2 \setminus \Delta) \to \Mod(T^2) \cong \SL_2\ZZ
\]
be the monodromy representation of $\pi$. By Moishezon's classification of elliptic fibrations (see \Cref{thm:moishezon} above), there are simple loops $\gamma_1,\ldots,\gamma_n$ about each critical value of $\Delta$ so that
\begin{align}
	\begin{split}
	\gamma_1\cdots \gamma_n &= \Id \\
	\phi_\pi(\gamma_{2i+1}) &= \begin{pmatrix} 1 & -1 \\ 0 & 1 \end{pmatrix}  \\
	\phi_\pi(\gamma_{2i}) &= \begin{pmatrix} 1 & 0 \\ 1 & 1 \end{pmatrix} .
	\end{split}
	\label{eq:standard-monodromy}
\end{align}
Furthermore, by \Cref{thm:moishezon-matsumuto}, a mapping class $f$ lies in $\Br(\pi)$ if and only if $\phi_\pi \circ f_\ast$ is conjugate to $\phi_\pi$ for any lift of $f$ to an automorphism of $\pi_1(S^2 \setminus \Delta)$ (i.e., via the Hurwitz action). We will denote that two representations $\phi,\phi'$ are conjugate by $\phi \simeq \phi'$. To simplify notation, let
\begin{align}
	r \coloneqq \sigma_1\cdots \sigma_{n-1} \in \Mod(S^2,\Delta) \label{eq:full-rotation}
\end{align}
be the $2\pi/n$ rotation of $(S^2,\Delta)$, where $\sigma_i$ is a half-twist supported on a neighborhood of the once-punctured disks on the interiors of $\gamma_i$ and $\gamma_{i+1}$. Note that, by the theorem of Murasugi mentioned in the introduction, the rotation $r$ represents the unique order $n$ conjugacy class in $\Mod(S^2,\Delta)$ \cite{murasugi-braid-torsion}. Before we proceed with the proof that part (a) of \Cref{thm:fin-order} and \Cref{thm:matrix-fin-order} are equivalent, we require an elementary fact about centralizers in spherical braid groups.
\begin{lemma}\label{lemma:centralizer}
	If $n \coloneqq \abs{\Delta} > 2$, then the centralizer $Z_{\Mod(S^2,\Delta)}(r)$ of $r = \sigma_1\cdots \sigma_{n-1}$ in $\Mod(S^2,\Delta)$ is the cyclic group $\langle r \rangle$.
\end{lemma}

\begin{proof}
	We deduce this from the analogous fact in the braid group. Let $B_n$ be the braid group on $n$ strands. Applying the capping homomorphism (see \cite[Section 3.6.2]{primer}), gives an exact sequence
	\begin{align*}
		1 \to \langle \widetilde{T}_n^2 \rangle \to B_n \to \Mod(S^2,\Delta) \to 1,
	\end{align*}
	with $\widetilde{T}_n$ the lift of $T_n$, defined by 
	\begin{align*}
		\widetilde{T}_n = (\widetilde{\sigma}_1\cdots \widetilde{\sigma}_{n-1})(\widetilde{\sigma}_1\cdots \widetilde{\sigma}_{n-2})\cdots (\widetilde{\sigma}_1\widetilde{\sigma}_2) \widetilde{\sigma}_1,
	\end{align*}
	where each $\widetilde{\sigma}_i$ is a lift of $\sigma_i$. Let $\widetilde{r}$ be the lift of $r$ which rotates the $n$ marked points around the center of the disk. Let $g \in \Mod(S^2,\Delta)$ centralize $r$, and choose a lift $\widetilde{g} \in B_n$. Then $\widetilde{g}\widetilde{r}\widetilde{g}^{-1} = \widetilde{r}T_n^{2k}$ for some $k \in \ZZ$. Furthermore, $\widetilde{r}^n = T_n^2$, and so this implies that $\widetilde{g}\widetilde{r}\widetilde{g}^{-1} = \widetilde{r}^{1+nk}$. Conjugation by $\widetilde{g}$ thus restricts to an automorphism of $\langle \widetilde{r} \rangle $, implying that $1+nk = \pm 1$. Since $k \in \ZZ$ and $n > 2$, the only solution to this equation is $k = 0$. Therefore $\widetilde{g}$ lies in the centralizer of $\widetilde{r}$. Proposition 3.3 of \cite{braid-centralizer} shows that $Z_{B_n}(\widetilde{r}) = \langle \widetilde{r}\rangle$, and so $\widetilde{g} \in \langle \widetilde{r} \rangle$. It follows that $g \in \langle r \rangle$ as claimed.
\end{proof}

With the above in hand we can now prove the equivalence of \Cref{thm:matrix-fin-order} and part (a) of \Cref{thm:fin-order}.
\begin{proof}[Proof that \Cref{thm:matrix-fin-order} implies part (a) of \Cref{thm:fin-order}]
	Let $r \in \Mod(S^2,\Delta)$ be the order $n$ rotation. By Murasugi's theorem (see the introduction), any order $n$ element of $\Br(\pi)$ is conjugate to $r$ by some $b \in \Mod(S^2,\Delta)$. Thus, it suffices to consider elements of the form $brb^{-1}$ for $b \in \Mod(S^2,\Delta)$. Choose loops $\gamma_1,\ldots,\gamma_n$ realizing the standard monodromy as in \eqref{eq:standard-monodromy}. A direct calculation shows that $\phi_\pi \circ r_\ast$ is conjugate to $\phi_\pi$ via the conjugating matrix $C_4 \coloneqq \begin{psmallmatrix} 0 & -1 \\ 1 & 0 \end{psmallmatrix}$ exchanging the simple closed curves corresponding to $\begin{psmallmatrix} 1 \\ 0 \end{psmallmatrix}$ and $\begin{psmallmatrix} 0 \\ 1 \end{psmallmatrix}$. By the aforementioned theorem of Moishezon (see \Cref{thm:moishezon-matsumuto}), the conjugate $br b^{-1}$ belongs to $\Br(\pi)$ if and only if $\phi_\pi \circ (brb^{-1})_\ast$ is conjugate to $\phi_\pi$ by some matrix in $\SL_2\ZZ$. Adopting the notation $\simeq$ for conjugacy between maps, we may rewrite this as
	\begin{align}
		\phi_\pi \circ b_\ast \circ r_\ast \circ b^{-1}_\ast &\simeq \phi_\pi \\
		\phi_\pi \circ b_\ast \circ r_\ast &\simeq \phi_\pi \circ b_\ast.\label{eq:br-conjugacy-equiv}
	\end{align}
	Let $X_i = \phi_\pi(b_\ast(\gamma_i))$. Then $X_1\cdots X_n = \Id$ since $b_\ast$ is an automorphism of $\pi_1$. Because $b_\ast$ is represented by a diffeomorphism, each $b_\ast(\gamma_i)$ is represented by a simple loop around a single puncture, and so Picard-Lefschetz theory implies that $X_i = T_v \in \SL_2\ZZ$ for some primitive vector $v \in \ZZ^2$. Direct calculation gives $r(\gamma_i) = \gamma_{i+1}$, so that
	\begin{align*}
		(\phi_\pi\circ b_\ast\circ r_\ast)(\gamma_i) = X_{i+1},
	\end{align*}
	and so the conjugacy \eqref{eq:br-conjugacy-equiv} implies that for some $C \in \SL_2\ZZ$ we have $X_{i+1} = CX_iC^{-1}$ for all $i$. In summary, $C,X_1,\ldots,X_n$ satisfy assumptions (1) to (3) in \Cref{thm:matrix-fin-order}:
	\begin{enumerate}
		\item $X_1 \in \SL_2\ZZ$ is given by symplectic transvection in the vanishing cycle, by Picard-Lefschetz theory.
		\item $X_{i+1} = CX_iC^{-1}$.
		\item $\prod_{i=1}^n X_i = \Id$.
	\end{enumerate}
	Thus, assuming that \Cref{thm:matrix-fin-order} holds, there exists some conjugating matrix $D \in \SL_2\ZZ$ so that either $X_i = DY_iD^{-1}$ for all $i$ or $X_i = DZ_iZ^{-1}$ for all $i$, where
	\begin{align*}
		(Y_1,Y_2,Y_3,\ldots) &\coloneqq \left(\begin{pmatrix} 1 & -1 \\ 0 & 1 \end{pmatrix}, \begin{pmatrix} 1 & 0 \\ 1 & 1 \end{pmatrix}, \begin{pmatrix} 1 & -1 \\ 0 & 1 \end{pmatrix},\ldots\right) \\
		(Z_1,Z_2,Z_3,Z_4,\ldots) &\coloneqq \left(\begin{pmatrix} 1 & -1 \\ 0 & 1 \end{pmatrix}, \begin{pmatrix} 2 & -1 \\ 1 & 0 \end{pmatrix}, \begin{pmatrix} 1 & 0 \\ 1 & 1 \end{pmatrix}, \begin{pmatrix} 1 & -1 \\ 0 & 1 \end{pmatrix}, \ldots\right).
	\end{align*}
	If $(X_1,\ldots,X_n)$ is simultaneously conjugate to $(Y_1,\ldots,Y_n)$, then $\phi_\pi \circ b_\ast \simeq \phi_\pi$ and so $b \in \Br(\pi)$. Therefore $br b^{-1} \in \Br(\pi)$ is an order $n$ element in the conjugacy class of $r$. If instead $(X_1,\ldots,X_n)$ is simultaneously conjugate to $(Z_1,\ldots,Z_n)$, then for the Garside twist $T_n$ defined by
	\begin{align*}
		T_n = (\sigma_1\sigma_2\cdots \sigma_{n-1})(\sigma_1\cdots \sigma_{n-2})\cdots(\sigma_1\sigma_2)\sigma_1
	\end{align*}
	we have that $\phi_\pi \circ b_\ast \simeq \phi_\pi \circ (T_n)_\ast$. Therefore $(T_nb^{-1})$ lies in $\Br(\pi)$ and so $br b^{-1}$ lies in the same conjugacy class as $r^{-1} = T_n r T_n^{-1}$. It is clear that $r^{-1}$ itself lies in $\Br(\pi)$, since $r$ does.

	Finally, it suffices to see that the conjugacy classes of $r$ and $r^{-1} = T_nr T_n^{-1}$ in $\Br(\pi)$ are distinct. Suppose not. Then there is some $b \in \Br(\pi)$ so that $br b^{-1} = T_n r T_n ^{-1}$. Then $T_n^{-1}b$ centralizes $r$, and so applying \Cref{lemma:centralizer}, we have that $T_n^{-1}b = r^k$ for some $k$. However, this would imply that $T_n = br^{-k} \in \Br(\pi)$, and we know that $(T_n \cdot \phi_\pi)(\gamma_i) = Z_i$ by direct calculation. The sequence $Z_i$ cannot be simultaneously conjugated to $Y_i$, and thus $T_n \not\in \Br(\pi)$, giving us a contradiction.
\end{proof}
The proof that part (a) of \Cref{thm:fin-order} implies \Cref{thm:matrix-fin-order} is similar to the above, and we omit it for brevity (and since we do not need it). The essential realization is that any matrices $X_i$ satisfying the assumptions of \Cref{thm:matrix-fin-order} induce a point in the Hurwitz orbit of $\phi_\pi$ which is fixed by the rotation $r$.

\section{Proving \texorpdfstring{\Cref{thm:matrix-fin-order}}{Theorem 1.2} via the Trace Polynomial}\label{sec:trace-poly}

\begin{proof}[Proof of \Cref{thm:matrix-fin-order}]
	Let $C, X_1,\ldots,X_n \in \SL_2\ZZ$ satisfy \Cref{item:prim,item:cyclic,item:prod-is-id} of \Cref{thm:matrix-fin-order}. For convenience, we recall \Cref{item:prim,item:cyclic,item:prod-is-id} here:
	\begin{enumerate}
		\item $X_1 = T_v$, where $T_v$ is the symplectic transvection $T_v(x) = x + (x,v)v$ for some primitive vector $v \in \ZZ^2$,
		\item $X_{i+1} = CX_iC^{-1}$ for $i \in \ZZ/n\ZZ$,
		\item $\prod_{i=1}^n X_i = \Id$.
	\end{enumerate}
	Note that $X_i = T_{C^{i-1}v}$, as $CT_wC^{-1} = T_{Cw}$ for any primitive vector $w \in \ZZ^2$. \\

	\emph{Step 1: The conjugating matrix $C$ has order $3$, $4$, or $6$.}  
 
	It is well known that $T_v = T_w$ if and only if $v = \pm w$, and so $T_v = T_{C^n v}$ implies $C^n v = \pm v$. As a consequence, each eigenvalue of $C$ is a root of unity. Since $C \in \SL_2\ZZ$, the eigenvalues of $C$ satisfy a monic degree two polynomial with integral coefficients, and so must be one of $1, -1, i, -i, \omega,$ or $-\omega$, where $\omega = e^{2\pi i/3}$ is a primitive 3rd root of unity. In the latter four cases, there are two distinct eigenvalues, and so $C^3, C^4, C^6 = \Id$ respectively in each case (by diagonalizing $C$). In the first case, where $C$ has eigenvalue $1$ with algebraic multiplicity $2$, this implies either $C = \Id$ or $Cv = v$, hence $X_i = T_v$ for all $i$, and so $\prod_{i=1}^n X_i = T_v^n \neq \Id$. Similarly, if $C$ has eigenvalue $-1$ then $X_i = T_v$, and the product cannot be the identity.\\

	Note further that if $C$ has order $4$ then $C^2 = -\Id$, and hence $X_i$ is periodic with period $2$. Otherwise, $X_i$ has period $3$. Note also that $(\SL_2 \ZZ)^{ab} = \ZZ/12\ZZ$ and $T_v$ maps to a generator of $\ZZ/12\ZZ$ under the abelianization map (see \cite[p. 123]{primer}). As a consequence, $n$ must be divisible by $12$. Since the period of $(X_1,\ldots,X_n)$ is $k \in \{2,3\}$, \Cref{item:prod-is-id} may be rewritten as 
	\begin{align*}
		(X_1\cdots X_k)^{n/k} = \Id.
	\end{align*}
	We define the \textit{trace polynomial} as
	\begin{align*}
		f(p,q) \coloneqq \tr(X_1\cdots X_k) = \tr(T_v \cdots T_{C^{k-1}}v) \in \ZZ[p,q],
	\end{align*}
	in the coordinates $v = \begin{psmallmatrix} p \\ q \end{psmallmatrix}$. Since $X_1 \cdots X_k$ is finite order, its trace has absolute value at most two, i.e., $\abs{f(p,q)} \leq 2$. \\

	\emph{Step 2: Up to a change of coordinates in $\SL_2\ZZ$, $f$ is given by an explicit polynomial of the norm 
	\[
		N(p+q\lambda) \coloneqq (p+q\lambda)(p+q\overline{\lambda})
	\]
	where $\lambda$ is a 3rd or 4th root of unity.}

	Because $C$ has order $3$,$4$, or $6$, it must be conjugate in $\SL_2\ZZ$ to one of the three matrices\footnote{These matrices correspond precisely to those regular polygons tiling the plane as well as the number rings $\ZZ[i]$ and $\ZZ[\omega]$} \begin{align*}
		C_3 = \begin{pmatrix} 
			0 & -1 \\
			1 & -1
		\end{pmatrix}, 
		C_4 = \begin{pmatrix} 
			0 & -1 \\ 
			1 & 0 
		\end{pmatrix}, 
		C_6 = \begin{pmatrix} 
			1 & -1 \\ 
			1 & 0
		\end{pmatrix},
	\end{align*}
	with orders $3$, $4$, and $6$ respectively, note that $C_3 = C_6^2$ (see \cite[p. 201]{primer}). Hence, by applying a global conjugacy, one may assume that $C$ is one of these three matrices. At this point, one can explicitly compute the polynomials in each case, which we will denote by $f_3,f_4,$ and $f_6$ respectively:
	\begin{align*}
		f_3(p,q) &= -(N_3+1)(N_3^2 + 2N_3 - 2) \\
		f_4(p,q) &= -2(N_4 + 1)(N_4 - 1)\\
		f_6(p,q) &= (N_3-1)(N_3^2 - 2N_3 - 2),
	\end{align*}
	where $N_j \coloneqq N(p+q\zeta_j)$ and $\zeta_j = e^{2\pi i/j}$ \\

	\textit{Step 3: Determining when $\abs{f_j(p,q)} \leq 2$ and $T_v \cdots T_{C^{k-1}v} = \Id$}.
	
	Take some $h \in \SL_2\ZZ$ so that $hCh^{-1} = C_j$ for $j \in \{3,4,6\}$. Then for $Y_i = hX_ih^{-1}$ we have $Y_1  = T_{hv}$, $Y_i = C_j Y_{i-1}C_j^{-1}$, and $\prod_{i=1}^n Y_i = \Id$. Thus we can assume without loss of generality that $C = C_j$, and instead classify the possible $v = \begin{psmallmatrix} p \\ q \end{psmallmatrix} $ such that $\abs{f_j(p,q)} \leq 2$. By direct computation, we find that 
	\begin{align*}
		\abs{\tr(X_1X_2X_3)} = \abs{f_3(p,q)} \leq 2 &\implies N_3 = N(p+q\omega) \leq 1 \tag{if $C = C_3$}\\
		\abs{\tr(X_1X_2)} = \abs{f_4(p,q)} \leq 2 &\implies N_4 = N(p+qi) < 2 \tag{if $C=C_4$} \\
		\abs{\tr(X_1X_2X_3)} = \abs{f_6(p,q)} \leq 2 &\implies N_3 = N(p+q\omega) \leq 3 \tag{if $C = C_6$},
	\end{align*}
	where $\omega = e^{2\pi i /3}$. These regions are the interiors of conics in the $(p,q)$ plane, as displayed in \Cref{fig:ordernconics-intro}. There are then a finite number of
	possibilities of $v \in \ZZ^2$ for each $C$. When $C = C_3$, there is one $C$-orbit of
	primitive vectors up to sign which satisfy $\abs{f_3(p,q)} \leq 2$:
	\begin{align*}
		v = \pm \begin{pmatrix} 
			1 \\ 
			0 
		\end{pmatrix} && 
		Cv = \pm \begin{pmatrix} 
			0 \\ 
			1 
		\end{pmatrix} &&
		C^2 v = \pm \begin{pmatrix} 
			-1 \\
			-1
		\end{pmatrix}.
	\end{align*}
	In this case, direct computation shows that $T_vT_{Cv}T_{C^2v}$ is a parabolic element of infinite order. When $C = C_4$, there is again one orbit up to sign
	\begin{align*}
		v = \pm \begin{pmatrix} 1 \\ 0 \end{pmatrix} && Cv = \pm \begin{pmatrix} 0 \\ 1 \end{pmatrix},
	\end{align*}
	in this case $T_vT_{Cv}$ has order $6$, and $T_v,T_{Cv}$ is precisely the pair $Y_1$ and $Y_2$ identified in the statement of the theorem,
	\begin{align*}
		Y_1 = \begin{pmatrix}
			1 & -1 \\
			0 & 1
		\end{pmatrix} && 
		Y_2 = \begin{pmatrix} 
			1 & 0 \\
			1 & 1 
		\end{pmatrix}.
	\end{align*}
	Finally, we turn to when $C = C_6$, there are then two $C$-orbits of primitive vectors up to sign where $N(p+q\omega) \leq 3$.
	\begin{align*}
		v &= \pm \begin{pmatrix} 
			1 \\
			0
		\end{pmatrix} && 
		C v = \pm
		\begin{pmatrix} 
			1 \\
			1
		\end{pmatrix} && 
		C^2 v = \pm
		\begin{pmatrix} 
			0 \\
			1
		\end{pmatrix},
		\\
		w &= \pm
		\begin{pmatrix} 
			2 \\
			1
		\end{pmatrix} &&
		Cw = \pm
		\begin{pmatrix} 
			1 \\
			2
		\end{pmatrix}  &&
		C^2w = \pm 
		\begin{pmatrix} 
			-1 \\
			1
		\end{pmatrix}.
	\end{align*}
	A simple calculation verifies that $T_vT_{Cv}T_{C^2v}$ has order $4$, and corresponds to the triple identified in the theorem statement
	\begin{align*}
		Z_1 = \begin{pmatrix} 1 & -1 \\0 & 1 \end{pmatrix} && Z_2 = \begin{pmatrix} 2 & -1 \\ 1 & 0 \end{pmatrix} && Z_3 \begin{pmatrix} 1 & 0 \\ 1 & 1 \end{pmatrix}.
	\end{align*}
	Furthermore, the product $T_wT_{Cw}T_{C^2w}$ is a parabolic element of infinite order. Because $(X_1\cdots X_k)$ must have finite order, these calculations imply that $X_1,\ldots,X_n$ must be equal to one of $Y_1,\ldots,Y_n$ or $Z_1,\ldots,Z_n$, depending on whether $C = C_4$ or $C = C_6$, completing the proof.
\end{proof}
Granted the equivalence of \Cref{thm:matrix-fin-order} and part (a) of \Cref{thm:fin-order} proved above, we have then completed the classification of order $n$ elements in $\Br(\pi)$ for $\pi : M \to S^2$ a genus one Lefschetz fibration with $n$ singular fibers.

\section{An Algebraic Approach to the Trace Polynomial}\label{sec:alg-trace-poly}

Before proving parts (b) and (c) of \Cref{thm:fin-order} and discussing the difficulties which arise from the approach above for orders less than $n-2$, we will describe an invariant theory approach to Step 2. This alternate description explains why one might expect the polynomials $f_3,f_4,f_6$ to be polynomials in the norms $N_4,N_3$ in more algebraic terms. To begin, we note the following proposition.
\begin{proposition}
	Let $C, X_1,\ldots,X_n \in \SL_2\ZZ$ satisfy \Cref{item:prim,item:cyclic,item:prod-is-id} of \Cref{thm:matrix-fin-order}, and let $(X_1,\ldots,X_n)$ have period $k$. Then
	\begin{align*}
		f(p,q) = \tr(T_vT_{Cv}\cdots T_{C^{k-1}v}) \in \ZZ[p,q],
	\end{align*}
	for $v = \begin{psmallmatrix} p \\ q \end{psmallmatrix}$, is invariant under the centralizer $Z_{\GL_2\ZZ}(C)$ of $C$ in $\GL_2 \ZZ$ acting on $\ZZ[p,q]$. Furthermore, $Z_{\GL_2\ZZ}(C) \cong D_{2k}$, the dihedral group with $2k$ elements.
\end{proposition}
Because $D_{2k}$ is a Coxeter group, standard methods such as the Chevalley-Shephard-Todd Theorem, which states that the invariant ring $\CC[p,q]^{D_{2k}}$ is a free polynomial algebra, apply over $\CC$. In this case, we are actually able to compute the full invariant ring over $\ZZ$, since there is an integral choice of generators for the invariant ring over $\CC$.
\begin{proposition}\label{prop:invariant-theory}
	Let $C \in \SL_2\ZZ$ have order $3, 4,$ or $6$ and let $H$ be the centralizer of $C$ in $\GL_2\ZZ$, then
	\begin{align*}
		\ZZ[p,q]^H \cong \begin{cases}
			\ZZ[N(x+y\omega), (x+y\omega)^6 + (x+y\overline{\omega})^6] & C \text{ has order } 3 \text{ or } 6 \\
			\ZZ[N(x+yi), x^2y^2] & C \text{ has order } 4
		\end{cases}
	\end{align*}
	where the isomorphism is by an $\SL_2\ZZ$ change of coordinates. 
\end{proposition}

\begin{proof}
	By applying a global conjugation, let $C$ be one of $C_3,C_4$ or $C_6$. The conjugating matrix $L$ so that $C = LC_jL^{-1}$ induces the coordinates $x,y$ given in the theorem. One verifies that the centralizer is generated by $C,-\Id, \begin{psmallmatrix} 0 & 1 \\ 1 & 0 \end{psmallmatrix}$ and is isomorphic to $D_{2k}$, where $k = 2$ if $C$ has order $4$ and $k = 3$ if $C$ has order $3$ or $6$. If $C = C_3$ then $C_6$ lies in the centralizer as well. A direct computation shows that the generators claimed above are in the invariant ring. The Jacobian criterion states that polynomials $g_1,\ldots,g_n \in \CC[x_1,\ldots,x_m]$ are algebraically independent provided that the differential $\d g_1 \wedge \cdots \d g_n$ is not identically zero \cite[\S 3.10]{humphreys-coxeter}. When $n = m$, we identify the wedge product with the determinant and write $J(g_1,\ldots,g_n) = \det \mathcal{J}(g_1,\ldots,g_n)$, where $\mathcal{J}(g_1,\ldots,g_n)$ is the matrix of partial derivatives. Applying the Jacobian criterion in this case verifies algebraic independence:
	\begin{align*}
		J(N(x+y\omega), (x+y\omega)^6 + (x-y\omega)^6) &= \det \begin{pmatrix} 
			2x - y & 6(2x-y)(x^4 - 2x^3y - 6x^2y^2 + 7xy^3 + 4) \\
			2y - x & 6(2y-x)(y^4 - 2y^3x - 6y^2x^2 + 7yx^3 + 4)
		\end{pmatrix} \\
		&= 6(2x-y)(2y-x)(y^4 - 9xy^3 + 9x^3y-x^4) \neq 0 \\
			J(N(x+yi),x^2y^2) = J(x^2+y^2,x^2y^2) &= 
			\det \begin{pmatrix} 
				2x & 2xy^2 \\
				2y & 2x^2y 
			\end{pmatrix} 
			= 2x^3y - 2xy^3 \neq 0.
	\end{align*}
	Let these two proposed generators be referred to as $F_1,F_2$. One may check that
	\begin{align*}
		\abs{D_{2k}} = 4k = \deg(F_1)\deg(F_2).
	\end{align*}
	Hence, these are generators for $\CC[x,y]^H$ (see \cite[p. 67]{humphreys-coxeter}). Because $F_1,F_2 \in \ZZ[x,y]$ are primitive polynomials, we in fact obtain that these are generators for $\ZZ[x,y]^H$ as desired.
\end{proof}
Note that the trace polynomial $f$ has total degree $4$ or $6$ depending on whether $(X_1,\ldots,X_k)$ has period $k = 2$ or $k = 3$. To show that $f(p,q)$ is a polynomial in the norm $N_4 = N(x+yi)$ or $N_3 = N(x+y\omega)$ respectively it thus suffices by \Cref{prop:invariant-theory} to determine the degree $2k$ homogeneous part of $f$. This computation must be carried out directly, and yields $-2N_4^2$ and $\pm N_3^3$ respectively.

%% file: orders-not-n.tex
\section{Orders \texorpdfstring{\bm{$n-1$}}{n-1} and \texorpdfstring{$\bm{n-2}$}{n-2}}\label{sec:orders-not-n}

Let $\tau \in \Mod(S^2,\Delta)$ be the order $n-1$ rotation and $\eta \in \Mod(S^2,\Delta)$ be the order $n-2$ rotation. For convenience, we modify $\eta$ by isotopy so that it fixes a neighborhood of the south pole, and choose a basepoint within this neighborhood. When dealing with $\tau$, we can choose the south pole itself as the basepoint, and let the north pole be one of the punctures. One can choose generators $\gamma_1',\ldots,\gamma_{n-1}', \de$ and $\gamma_1'',\ldots,\gamma''_{n-2}, \nu_1, \nu_2$ for $\pi_1(S^2 \setminus \Delta)$ so that
\begin{align*}
	\tau(\gamma_i') &= \gamma_{i+1}' \tag{for $i \in \ZZ/(n-1)\ZZ$} \\
	\tau(\de) &= \gamma_1^{-1}\de\gamma_1 \\
	\eta(\gamma_i'') &= \gamma_{i+1}'' \tag{for $i \in \ZZ/(n-2)\ZZ$} \\
	\eta(\nu_1) &= \nu_1 \\
	\eta(\nu_2) &= \gamma_1^{-1}\nu_2\gamma_1,
\end{align*}
where $\gamma_i'$ and $\gamma_i''$ are loops about the equatorial punctures, $\de$ and $\nu_1$ are loops about the north pole, and $\nu_2$ is a loop about the south pole. Note that, choosing $\gamma_i', \gamma_i'',\de,\nu_1,$ and $\nu_2$ appropriately gives
\begin{align*}
	\de\gamma_1'\cdots \gamma_{n-1}' &= 1 && \nu_1 \gamma_1'' \cdots \gamma_{n-2}'' \nu_2 = 1.
\end{align*}
Similarly to \Cref{sec:equiv}, to prove parts (b) and (c) of \Cref{thm:fin-order} it suffices to show that there are no monodromy representations which are fixed up to conjugacy by $\tau$ and $\eta$ respectively. By evaluating such a representation $\rho_\pi : \pi_1(S^2 \setminus \Delta) \to \SL_2\ZZ$ at our choice of generators, we reduce to showing that there are no monodromy factorizations fixed by $\tau$ or $\eta$. Thus, parts (b) and (c) of \Cref{thm:fin-order} are implied by the following statements.
\begin{theorem}\label{thm:matrix-fin-order-n-minus-one}
	Let $n \geq 1$, then there are no tuples $(X_1,\ldots,X_{n-1},L) \in (\SL_2\ZZ)^n$ of matrices satisfying the following for some $C \in \SL_2\ZZ$:
	\begin{enumerate}[(1')]
		\item $X_1$ and $L$ are symplectic transvections,
		\item $L\prod_{i=1}^{n-1}X_i = \Id$,
		\item $CX_iC^{-1} = X_{i+1}$ and $CLC^{-1} = X_1^{-1}LX_1$.
	\end{enumerate}
\end{theorem}

\begin{proposition}\label{prop:matrix-fin-order-n-minus-two}
	Let $n \geq 1$, then there are no tuples $(X_1,\ldots,X_{n-2},L_1,L_2) \in (\SL_2\ZZ)^n$ satisfying the following for some $C \in \SL_2\ZZ$:
	\begin{enumerate}[(1'')]
		\item $X_1,L_1,L_2$ are symplectic transvections,
		\item $L_1\left(\prod_{i=1}^{n-2}X_i\right)L_2 = \Id$,
		\item $CX_iC^{-1} = X_{i+1}$, $CL_1C^{-1} = X_1^{-1}L_1X_1$, and $CL_2C^{-1} = L_2$.
	\end{enumerate}
\end{proposition}

We now show that \Cref{thm:matrix-fin-order-n-minus-one} and \Cref{prop:matrix-fin-order-n-minus-two} imply parts (b) and (c) of \Cref{thm:fin-order} respectively.
\begin{proof}[Proof that \Cref{thm:matrix-fin-order-n-minus-one} implies part (b) of \Cref{thm:fin-order}]j
	Let $\pi : M \to S^2$ be an elliptic fibration, and let $\tau \in \Mod(S^2,\Delta)$ be the order $n-1$ rotation. By Murasugi's theorem (see the introduction), $\tau$ represents the unique order $n-1$ conjugacy class in $\Mod(S^2,\Delta)$ \cite{murasugi-braid-torsion}. Hence any order $n-1$ element of $\Br(\pi)$ is represented by $b\tau b^{-1} \in \Br(\pi)$ for some $b \in \Mod(S^2,\Delta)$. Let 
	\[
		\rho_\pi : \pi_1(S^2 \setminus \Delta) \to \SL_2\ZZ
	\]
	be the monodromy representation of $\pi$.  By Moishezon's theorem (see \Cref{thm:moishezon-matsumuto} above), we know that $b\tau b^{-1} \in \Br(\pi)$ if and only if $\rho_\pi \circ (b\tau b^{-1})_\ast$ is conjugate to $\rho_\pi$. We choose generators $\gamma_1',\ldots,\gamma_{n-1}', \de$ so that
	\begin{align*}
		\tau(\gamma_i') &= \gamma_{i+1}' \tag{for $i \in \ZZ/(n-1)\ZZ$} \\
		\tau(\de) &= \gamma_1^{-1}\de\gamma_1 \\
	\end{align*}
	where $\gamma_i'$ are loops about the equatorial punctures and $\de$ is a loop about the marked fixed point of $\tau$. Note that $\de \gamma_1'\cdots \gamma_{n-1}' = \Id$.  Therefore, letting $X_i = \rho_\pi(b_\ast(\gamma_i'))$ and $L = \rho_\pi(b_\ast(\de))$ we conclude that there is some matrix $C \in \SL_2\ZZ$ so that
	\begin{align*}
		CX_iC^{-1} &= C\rho_\pi(b_\ast(\gamma_i'))C^{-1} = \rho_\pi(b_\ast(\tau_\ast(\gamma_i'))) = \rho_\pi(b_\ast(\gamma_{i+1}')) = X_{i+1} \\
		CLC^{-1} &= C\rho_\pi(b_\ast(\de))C^{-1} = \rho_\pi(b_\ast(\tau_\ast(\de))) = \rho_\pi(b_\ast(\gamma_1^{-1}\de\gamma_1)) = X_1^{-1}LX_1.
	\end{align*}
	Furthermore, because $\rho_\pi$ and $b_\ast$ are group homomorphisms, $LX_1\cdots X_n = \Id$. By Picard-Lefschtz theory, we also know that $X_i,L$ are symplectic transvections. Therefore, $(X_1,\ldots,X_{n-1},L) \in (\SL_2\ZZ)^n$ is a tuple satisfying the conditions of \Cref{thm:matrix-fin-order-n-minus-one}. No such tuple exists, and so there can be no such $b \in \Mod(S^2,\Delta)$ so that $b\tau b^{-1} \in \Br(\pi)$.
\end{proof}
The proof that \Cref{prop:matrix-fin-order-n-minus-two} implies part (c) of \Cref{thm:fin-order} is identical, and so we omit it for brevity. With this motivation, we prove \Cref{thm:matrix-fin-order-n-minus-one} and \Cref{prop:matrix-fin-order-n-minus-two}.
\begin{proof}[Proof of \Cref{thm:matrix-fin-order-n-minus-one}]
	Let $(X_1,\ldots,X_{n-1},L)$ satisfy the conditions of the theorem and let $X_1 = T_v$ and $L = T_w$ for primitive vectors $v,w \in \ZZ^2$. As before, $n = 12d$ for some $d \geq 1$ since $(\SL_2\ZZ)^{ab} = \ZZ/12\ZZ$. Note that $C^{n-1}$ fixes $X_1$ by conjugation. Hence $C^{n-1}(v) = \pm v$. As argued in Step 1 in \Cref{sec:trace-poly}, either $C \cdot v = \pm v$ or $C$ has order $3,4$, or $6$. If $C \cdot v = \pm v$, then $X_1 = X_i$ for all $i$. Therefore $L = X_1^{-n+1}$. However, since $n \geq 12$, this is impossible, since the only symplectic transvection which is a power of $X_1$ is $X_1$ itself.

	If $C^2 = -\Id$, then
	\begin{align*}
		T_w = C^2T_wC^{-2} = X_2^{-1}X_1^{-1}T_wX_1X_2 = T_{(X_1X_2)^{-1}w}.
	\end{align*}
	and so $X_1X_2(w) = \pm w$. Therefore $X_1X_2 = \pm L^k$ for some $k \in \ZZ$, and so
	\begin{align*}
		L(X_1X_2)^{6d}X_2^{-1} = \Id,
	\end{align*}
	but then $\pm L^{6dk+1} = X_2$. Because $X_2$ and $L$ are both symplectic transvections, this is impossible unless $6dk+1 = 1$. However, this implies that $k=0$, and so $X_1X_2 = \Id$. Applying the abelianization map, this would imply that $2 = 0$ in $(\SL_2\ZZ)^{ab} \cong \ZZ/12\ZZ$, a contradiction.

	Otherwise, $C^3 = \pm \Id$ and
	\begin{align*}
		T = C^3T_wC^{-3} = X_3^{-1}X_2^{-1}X_1^{-1}T_wX_1X_2X_3 = T_{(X_1X_2X_3)^{-1}w}
	\end{align*}
	and so $X_1X_2X_3(w) = \pm w$. As before, $X_1X_2X_3 = \pm L^k$ for some $k \in \ZZ$ and so
	\begin{align*}
		L(X_1X_2X_3)^{4d}X_3^{-1} = \Id,
	\end{align*}
	which implies $\pm L^{4dk+1} = X_3$. As before, this implies that $3=0$ in $\ZZ/12\ZZ$.
\end{proof}

\begin{proof}[Proof of \Cref{prop:matrix-fin-order-n-minus-two}]
	Let $(X_1,\ldots,X_{n-2},L_1,L_2)$ satisfy the conditions of the proposition. Conjugating Item (2'') by the given $C \in \SL_2\ZZ$ gives
	\begin{align*}
		X_1^{-1}L_1 X_1 \cdots X_{n-2} X_1 L_2 &= \Id \\
		X_1^{-1} L_2^{-1}X_1L_2 &= \Id,
	\end{align*}
	and so $L_2^{-1}X_1L_2 = X_1$. Let $X_1 = T_v$, where $v \in \ZZ^2$ is some primitive vector. Then $L_2 \cdot v = \pm v$, and since $L_2$ is a symplectic transvection $L_2 = X_1$. Because $CL_2C^{-1} = L_2$, we also have that $X_i = X_1$ for all $i$. Therefore
	\begin{align*}
		L_1X_1^{n-1} = L_1X_1\cdots X_{n-2} L_2 = \Id,
	\end{align*}
	and so $L_1 = X_1^{1-n}$. Since $L_1$ and $X_1$ are both symplectic transvections, this is only possible if $1-n = 1$. But then $n = 0$, and so the only such tuple is the empty tuple.
\end{proof}
Hence, we have proved parts (b) and (c) of \Cref{thm:fin-order}.

%% file: tikz/farey-complex.tex
    \begin{tikzpicture}
		\begin{scope}[scale=1/4]
        \draw [thick] (-8,0) -- (8,0);
        \draw [thick] (0,0) circle (8);
        \foreach \i in {0,1,2,3} {%
            \draw [thick] (90*\i:8) arc (270+90*\i:180+90*\i:8);}
        \foreach \i in {0,1,...,7} {%
            \draw  (45*\i:8) arc (270+45*\i:135+45*\i:3.3);}
        \foreach \i in {0,1,...,15} {%
            \draw [thin] (22.5*\i:8) arc (270+22.5*\i:112.5+22.5*\i:1.6);}
        \foreach \i in {0,1,...,31} {%
            \draw [ultra thin] (11.25*\i:8) arc (270+11.25*\i:101.25+11.25*\i:0.8);}
        %\foreach \i in {0,1,...,63} {%
            %\draw [ultra thin] (5.625*\i:8) arc (270+5.625*\i:95.625+5.625*\i:0.4);}
		\end{scope}
		\node[circle,draw=red, fill=red, inner sep=0, minimum size=3.5pt, label=left:{$\color{red}\alpha$}] (a) at (-2,0) {};
		\node[circle,draw=blue, fill=blue, inner sep=0, minimum size=3.5pt, label=right:{$\color{blue}\beta$}] (b) at (2,0) {};
		\node[circle,draw=teal, fill=teal, inner sep=0, minimum size=3.5pt, label=above:{$\color{teal}\delta$}] (c) at (0,2) {};
		\node[circle,draw=violet, fill=violet, inner sep=0, minimum size=3.5pt, label=below:{$\color{violet}\gamma$}] (c) at (0,-2) {};
		\node[circle,draw=orange,fill=orange,inner sep=0, minimum size=3.5pt, label=below right:{\color{orange}$i$}] (i) at (0,0) {};
		\node[circle,draw=magenta,fill=magenta,inner sep=0, minimum size=3.5pt, label=right:{\color{magenta}$\omega$}] (w) at (0,0.5) {};
    \end{tikzpicture}